# A new proof of Menelaus's Theorem of Hyperbolic Quadrilaterals in the Poincaré Model of Hyperbolic Geometry

Florentin Smarandache [1] and Cătălin Barbu [2]

[1] Department of Mathematics, University of New Mexico, Gallup, NM 87301, USA.
E-mail: smarand@unm.edu

[2] Vasile Alecsandri College - Bacău, str. Vasile Alecsandri, nr. 37, cod 600011, Romania, E-mail: kafka_mate@yahoo.com

In this study, we present a proof of the Menelaus theorem for quadrilaterals in hyperbolic geometry, and a proof for the transversal theorem for triangles

**2000 Mathematics Subject Classification:** 51K05, 51M10

**Key words**: hyperbolic geometry, hyperbolic quadrilateral, Menelaus theorem, the transversal theorem, gyrovector

## 1 Introduction

Hyperbolic geometry appeared in the first half of the $19^{th}$ century as an attempt to understand Euclid's axiomatic basis for geometry. It is also known as a type of non-Euclidean geometry, being in many respects similar to Euclidean geometry. Hyperbolic geometry includes such concepts as: distance, angle and both of them have many theorems in common.There are known many main models for hyperbolic geometry, such as: Poincaré disc model, Poincaré half-plane, Klein model, Einstein relativistic velocity model, etc. The hyperbolic geometry is a non-euclidian geometry. Menelaus of Alexandria was a Greek mathematician and astronomer, the first to recognize geodesics on a curved surface as natural analogs of straight lines. Here, in this study, we give hyperbolic version of Menelaus theorem for quadrilaterals. The well-known Menelaus theorem states that if $l$ is a line not through any vertex of a triangle $ABC$ such that $l$ meets $BC$ in $D$, $CA$ in $E$, and $AB$ in $F$, then $\frac{DB}{DC} \cdot \frac{EC}{EA} \cdot \frac{FA}{FB} = 1$ [1]. We use in this study of Poincaré disc model.

We begin with the recall of some basic geometric notions and properties in the Poincaré disc. Let $D$ denote the unit disc in the complex $z$ - plane, i.e.
$$D = \{z \in \mathbb{C} : |z| < 1\}.$$
The most general Möbius transformation of $D$ is
$$z \to e^{i\theta} \frac{z_0 + z}{1 + \overline{z_0} z} = e^{i\theta}(z_0 \oplus z),$$
which induces the Möbius addition $\oplus$ in $D$, allowing the Möbius transformation of the disc to be viewed as a Möbius left gyro-translation

$$z \to z_0 \oplus z = \frac{z_0 + z}{1 + \overline{z_0} z}$$

followed by a rotation. Here $\theta \in \mathbb{R}$ is a real number, $z, z_0 \in D$, and $\overline{z_0}$ is the complex conjugate of $z_0$. Let $Aut(D, \oplus)$ be the automorphism group of the grupoid $(D, \oplus)$. If we define

$$gyr : D \times D \to Aut(D, \oplus), gyr[a,b] = \frac{a \oplus b}{b \oplus a} = \frac{1 + a\overline{b}}{1 + \overline{a}b},$$

then is true gyro-commutative law

$$a \oplus b = gyr[a,b](b \oplus a).$$

A gyro-vector space $(G, \oplus, \otimes)$ is a gyro-commutative gyro-group $(G, \oplus)$ that obeys the following axioms:

 (1) $gyr[\mathbf{u}, \mathbf{v}]\mathbf{a} \cdot gyr[\mathbf{u}, \mathbf{v}]\mathbf{b} = \mathbf{a} \cdot \mathbf{b}$ for all points $\mathbf{a}, \mathbf{b}, \mathbf{u}, \mathbf{v} \in G$.

 (2) $G$ admits a scalar multiplication, $\otimes$, possessing the following properties. For all real numbers $r, r_1, r_2 \in \mathbb{R}$ and all points $\mathbf{a} \in G$:

 (G1) $1 \otimes \mathbf{a} = \mathbf{a}$
 (G2) $(r_1 + r_2) \otimes \mathbf{a} = r_1 \otimes \mathbf{a} \oplus r_2 \otimes \mathbf{a}$
 (G3) $(r_1 r_2) \otimes \mathbf{a} = r_1 \otimes (r_2 \otimes \mathbf{a})$
 (G4) $\dfrac{|r| \otimes \mathbf{a}}{\|r \otimes \mathbf{a}\|} = \dfrac{\mathbf{a}}{\|\mathbf{a}\|}$
 (G5) $gyr[\mathbf{u}, \mathbf{v}](r \otimes \mathbf{a}) = r \otimes gyr[\mathbf{u}, \mathbf{v}]\mathbf{a}$
 (G6) $gyr[r_1 \otimes \mathbf{v}, r_1 \otimes \mathbf{v}] = 1$

 (3) Real vector space structure $(\|G\|, \oplus, \otimes)$ for the set $\|G\|$ of one-dimensional "vectors"

$$\|G\| = \{\pm \|\mathbf{a}\| : \mathbf{a} \in G\} \subset \mathbb{R}$$

with vector addition $\oplus$ and scalar multiplication $\otimes$, such that for all $r \in \mathbb{R}$ and $\mathbf{a}, \mathbf{b} \in G$,

 (G7) $\|r \otimes \mathbf{a}\| = |r| \otimes \|\mathbf{a}\|$
 (G8) $\|\mathbf{a} \oplus \mathbf{b}\| \leq \|\mathbf{a}\| \oplus \|\mathbf{b}\|$

**Definition 1** *The hyperbolic distance function in $D$ is defined by the equation*

$$d(a,b) = |a \:!\: b| = \left|\frac{a-b}{1-\overline{a}b}\right|.$$

Here, $a \:!\: b = a \oplus (-b)$, for $a, b \in D$.

For further details we refer to the recent book of A.Ungar [2].

**Theorem 2** ( *The Menelaus's Theorem for Hyperbolic Gyrotriangle* ). *Let $ABC$ be a gyrotriangle in a Möbius gyrovector space $(V_s, \oplus, \otimes)$ with vertices $A, B, C \in V_s$, sides $\mathbf{a}, \mathbf{b}, \mathbf{c} \in V_s$, and side gyrolengths $a, b, c \in (-s, s)$, $\mathbf{a} = !\ B \oplus C$, $\mathbf{b} = !\ C \oplus A$, $\mathbf{c} = !\ A \oplus B$, $a = \|\mathbf{a}\|$, $b = \|\mathbf{b}\|$, $c = \|\mathbf{c}\|$, and with gyroangles $\alpha, \beta,$ and $\gamma$ at the vertices $A, B,$ and $C$. If $l$ is a gyroline not through any vertex of an gyrotriangle $ABC$ such that $l$ meets $BC$ in $D$, $CA$ in $E$, and $AB$ in $F$, then*

$$\frac{(AF)_\gamma}{(BF)_\gamma} \cdot \frac{(BD)_\gamma}{(CD)_\gamma} \cdot \frac{(CE)_\gamma}{(AE)_\gamma} = 1.$$

where $v_\gamma = \dfrac{v}{1 - \dfrac{v^2}{s^2}}$ [3].

## 2  Main results

In this section, we prove Menelaus's theorem for hyperbolic quadrilateral.

**Theorem 3** ( *The Menelaus's Theorem for Gyroquadrilateral* ). *If $l$ is a gyroline not through any vertex of a gyroquadrilateral $ABCD$ such that $l$ meets $AB$ in $X$, $BC$ in $Y$, $CD$ in $Z$, and $DA$ in $W$, then*

$$\frac{(AX)_\gamma}{(BX)_\gamma} \cdot \frac{(BY)_\gamma}{(CY)_\gamma} \cdot \frac{(CZ)_\gamma}{(DZ)_\gamma} \cdot \frac{(DW)_\gamma}{(AW)_\gamma} = 1. \tag{1}$$

*Proof.* Let $T$ be the intersection point of the gyroline $DB$ and the gyroline $XYZ$ (See Figure 1). If we use a theorem 2 in the gyrotriangles $ABD$ and $BCD$ respectively, then

$$\frac{(AX)_\gamma}{(BX)_\gamma} \cdot \frac{(BT)_\gamma}{(DT)_\gamma} \cdot \frac{(DW)_\gamma}{(AW)_\gamma} = 1 \tag{2}$$

and

$$\frac{(DT)_\gamma}{(BT)_\gamma} \cdot \frac{(CZ)_\gamma}{(DZ)_\gamma} \cdot \frac{(BY)_\gamma}{(CY)_\gamma} = 1. \tag{3}$$

Multiplying relations (2) and (3) member with member, we obtain

$$\frac{(AX)_\gamma}{(BX)_\gamma} \cdot \frac{(BY)_\gamma}{(CY)_\gamma} \cdot \frac{(CZ)_\gamma}{(DZ)_\gamma} \cdot \frac{(DW)_\gamma}{(AW)_\gamma} = 1.$$

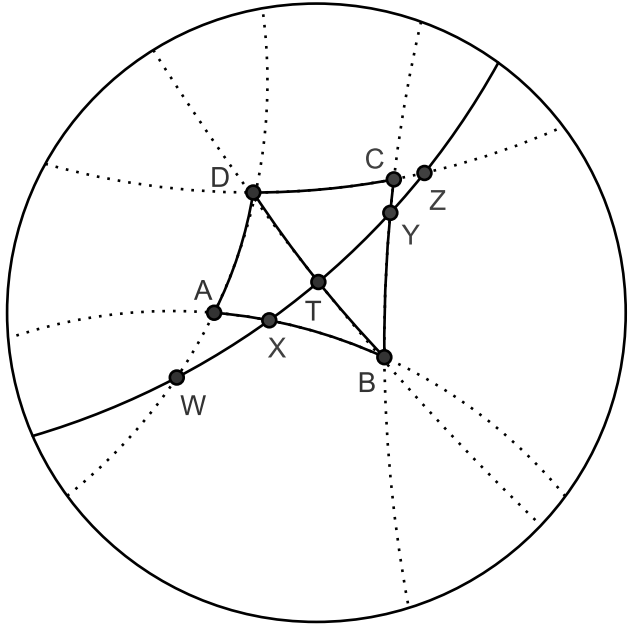

Figure 1

Naturally, one may wonder whether the converse of Menelaus theorem for hyperbolic quadrilateral exists. Indeed, a partially converse theorem does exist as we show in the following theorem.

**Theorem 4** ( *Converse of Menelaus's Theorem for Gyroquadrilateral* ). *Let $ABCD$ be a gyroquadrilateral. Let the points $X, Y, Z,$ and $W$ be located on the gyrolines $AB, BC, CD,$ and $DA$ respectively. If three of four gyropoints $X, Y, Z, W$ are collinear and*

$$\frac{(AX)_\gamma}{(BX)_\gamma} \cdot \frac{(BY)_\gamma}{(CY)_\gamma} \cdot \frac{(CZ)_\gamma}{(DZ)_\gamma} \cdot \frac{(DW)_\gamma}{(AW)_\gamma} = 1,$$

then all four gyropoints are collinear.

*Proof.* Let the points $W, X, Z$ are collinear, and gyroline $WXZ$ cuts gyroline $BC$, at $Y'$ say. Using the already proven equality (1), we obtain

$$\frac{(AX)_\gamma}{(BX)_\gamma} \cdot \frac{(BY')_\gamma}{(CY')_\gamma} \cdot \frac{(CZ)_\gamma}{(DZ)_\gamma} \cdot \frac{(DW)_\gamma}{(AW)_\gamma} = 1,$$

then we get

$$\frac{(BY)_\gamma}{(CY)_\gamma} = \frac{(BY')_\gamma}{(CY')_\gamma}. \tag{4}$$

This equation holds for $Y = Y'$. Indeed, if we take $x := \left| ! \ B \oplus Y' \right|$ and $b := \left| ! \ B \oplus C \right|$, then

we get $b \ominus x = |\ominus Y' \oplus C|$. For $x \in (-1,1)$ define

$$f(x) = \frac{x}{1-x^2} : \frac{b \ominus x}{1-(b \ominus x)^2}. \quad (5)$$

Because $b \ominus x = \dfrac{b-x}{1-bx}$, then $f(x) = \dfrac{x(1-b^2)}{(b-x)(1-bx)}$. Since the following equality holds

$$f(x) - f(y) = \frac{b(1-b^2)(1-xy)}{(b-x)(1-bx)(b-y)(1-by)}(x-y), \quad (6)$$

we get $f(x)$ is an injective function and this implies $Y = Y'$, so $W, X, Z,$ and $Y$ are collinear.

We have thus obtained in (1) the following

**Theorem 5** (*Transversal theorem for gyrotriangles*). *Let $D$ be on gyroside $BC$, and $l$ is a gyroline not through any vertex of a gyrotriangle $ABC$ such that $l$ meets $AB$ in $M$, $AC$ in $N$, and $AD$ in $P$, then*

$$\frac{(BD)_\gamma}{(CD)_\gamma} \cdot \frac{(CA)_\gamma}{(NA)_\gamma} \cdot \frac{(NP)_\gamma}{(MP)_\gamma} \cdot \frac{(MA)_\gamma}{(BA)_\gamma} = 1. \quad (7)$$

*Proof.* If we use a theorem 2 for gyroquadrilateral $BCNM$ and collinear gyropoints $D, A, P,$ and $A$ (See Figure 2), we obtain the conclusion.

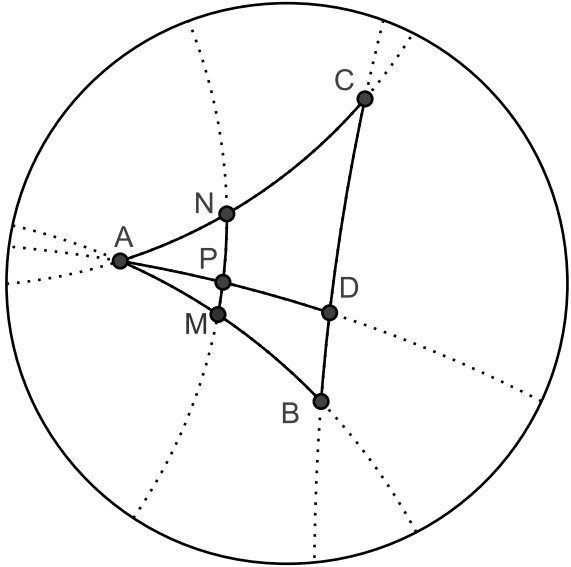

Figure 2

The Einstein relativistic velocity model is another model of hyperbolic geometry. Many of the theorems of Euclidean geometry are relatively similar form in the Poincaré model, Menelaus's theorem for hyperbolic gyroquadrilateral and the transversal theorem for gyrotriangle are an examples in this respect. In the Euclidean limit of large $s$, $s \to \infty$, gamma factor $v_\gamma$ reduces to $v$, so that the gyroinequalities (1) and (7) reduces to the

$$\frac{AX}{BX} \cdot \frac{BY}{CY} \cdot \frac{CZ}{DZ} \cdot \frac{DW}{AW} = 1$$

and

$$\frac{BD}{CD} \cdot \frac{CA}{NA} \cdot \frac{NP}{MP} \cdot \frac{MA}{BA} = 1,$$

in Euclidean geometry. We observe that the previous equalities are identical with the equalities of theorems of euclidian geometry.